\documentclass[12pt]{article}

\usepackage{amsmath}
\usepackage{amsfonts}
\usepackage{amssymb}

\newtheorem{Teo}{\textsc{Theorem}}
\newtheorem{Lemma}{{\textsc Lemma}}

\newtheorem{Cor}{\textsc{Corollary}}

\newcommand{\QED}{\hfill $\square$\bigskip}
\newcommand{\qed}{\hfill $\square$\bigskip}
\newcommand{\pf}{{\noindent\bf Proof.\ \ }}

\newcommand{\semidir}{\kern-2pt>\kern-6pt\triangleleft}
\newcommand{\Aut}{\mbox{\rm Aut}}
\newcommand{\Inn}{\mbox{\rm Inn}}
\newcommand{\N}{{\Bbb N}}

\newcommand{\la}{{\langle}}
\newcommand{\ra}{{\rangle}}
\newcommand{\g}{{\gamma}}%
\newcommand{\TI}{{${\bf T}\kern-10pt\sim$}}

\newcommand{\n}{\lhd}
\newcommand{\sn}{{\rm \kern1pt sn \kern1pt}}
\newcommand{\an}{{\rm \kern2pt an \kern1pt}}
\newcommand{\cf}{{\rm \kern2pt cf \kern1pt}}
\newcommand{\f}{{\rm \kern2pt f \kern1pt }}
\newcommand{\ff}{{\kern1pt \rm \bar{f}\kern1pt }}
\newcommand{\C}{\rm \kern2pt c \kern1pt}
\newcommand{\cn}{\rm \kern2pt cn\kern1pt}
\newcommand{\ssim}{\rm \kern2pt \approx}

\begin{document}



\vskip-1cm
\title{\vskip-1cm
Variants of theorems of  Baer and Hall\\ on finite-by-hypercentral groups
}
\vskip-1cm
\author{Carlo Casolo\ - Ulderico Dardano\  -\
Silvana Rinauro
 }

\date{\sl \small dedicated to the memory of Guido Zappa}

\vskip-2cm
\maketitle
 
\noindent\textbf{Abstract} {\small \ We show that if a group $G$ has a finite normal subgroup $L$ such that $G/L$ is hypercentral,  then the index of the hypercenter of $G$ is bounded by a function of the order of $L$. This completes recent results generalizing classical theorems by R. Baer and P. Hall. Then we apply our results to groups of automorphisms of a group $G$ acting in a restricted way on an ascending normal series of $G$.}\\

\section{Introduction}

A classical theorem by R. Baer states that, if the $m$-th term $Z_m(G)$ of the upper central series a group $G$ has finite index $t$ in $G$ for some positive integer $m$, then there is a finite normal subgroup $L$ of $G$ such that $G/L$ is nilpotent of class at most $m$, that is  $G/L=Z_m(G/L)$ (see 14.5.1 in \cite{RB}, which shall be the reference for undefined notation). Recently, in \cite{KS} it has been shown that there is such an $L$ with finite order $d$ bounded by a function of $t$ and $m$. 

In the opposite direction, P. Hall showed that, if there is a normal subgroup $L$ with finite order $d$ such that $G/L$ is nilpotent of class at most $m$, then $G/Z_{2m}(G)$ has finite order  bounded by a function of $d$ and $m$  (see \cite{RB}, page 118). 

Recently, in \cite{DGMS} it has been shown that the hypercenter  of $G$ has finite index $t$ if and only if  
there is a finite normal subgroup $L$ with order $d$  such that $G/L$ is { hypercentral}, that is coincides with its hypercenter.  
Recall that the {\it hypercenter} of a group $G$ is 
the last term of the upper central series of G (see details below). 
Then in \cite{KOS} it has been shown that $d$ may be bounded by a function of $t$, namely $t^{(1+\log_2 t)/2}$.
 Here we complete the picture by showing that $t$ in turn may be bounded by a function of  $d$. 

\begin{Teo}\label{TeoremaCasolo}
If a group $G$ has a finite normal subgroup $L$ such that $G/L$ is hypercentral,  then  
the hypercenter of $G$ has index bounded by 
$|\Aut(L)|\cdot |Z(L)|$. 
\end{Teo}

\begin{Cor}\label{CorTeoremaCasolo}
If a group $G$ has a finite normal subgroup $L$ such that $G/L$ is nilpotent of class $m$,  then  $|G/Z_{2m}(G)|$ is bounded by a function of $d:=|L|$. 
\end{Cor}

There are many generalization and variants of Baer and Hall theorems. By applying  Theorem \ref{TeoremaCasolo} above, we improve the results in \cite{DKP} which are concerned with possibly non-inner automorphisms. 

Before stating our Theorem 2 we recall some definition.  As usual, we say that the group $A$ \emph{acts} on a group $G$ if and only if   there is a homomorphism $\tilde{}: A\to \Aut(G)$ (called \emph{action}). 
We will regard both $G$ and $\tilde A$ as subgroups of the holomorph group $G\rtimes \Aut(G)$ of $G$. In particular, we will denote by a bar \ $\bar{}$\ \ the action of a group $G$ on itself by conjugation, that is the natural $\Aut(G)$-homomorphism $G\to \bar G\le \Aut(G)$. 
 If an action is such that 
 its image $\tilde {A}$ is normalized by $\bar G=\Inn(G)$, we define by recursion an ascending $G$-series  $Z_\alpha(G,A)$ (with $\alpha$ ordinal number) by $Z_0(G,A):=1$, $Z_{\alpha+1}(G,A)/Z_{\alpha}(G,A):=C_{G/Z_{\alpha}(G,A)}(A)$ and $Z_\lambda(G,A):=\cup_{\alpha<\lambda}Z_{\alpha}(G,A)$ when $\lambda$ is a limit ordinal.  We call $Z_\alpha(G,A)$ the $\alpha$th $A$-{\it center} of $G$. Recall that an ascending $G$-series is a well ordered (by inclusion) set of normal  subgroups of $G$. Clearly { the series $Z_\alpha(G,A)$ is stabilized by $A$}, in the sense that $A$ acts trivially on the factors beetwen consecutive terms. The last term $Z_\infty(G,A)$ of this series is called $A$-{\it hypercenter} of $G$.

We say that {\it $G$ is $A$-hypercentral with (ordinal) type at most $\alpha$} if and only if  
$G=Z_\alpha(G,A)$. Clearly $Z_\alpha(G):=Z_\alpha(G,\bar G)$ is the usual $\alpha$th  center of $G$ and if
 $G=Z_\alpha(G)$, then $G$ is hypercentral of type at most $\alpha$.

\medskip

Now we are in a position to state our second result, which consists in two parts that  refer to theorems of Baer and Hall, respectively. In fact, if $A=\Inn(G)$, then part (B) reduces to Theorem B in \cite{KOS} and part (H) to Theorem \ref{TeoremaCasolo} above. 

\begin{Teo}\label{generalizeBH} 
Let  $G$ be a group and $A$ be a subgroup of $\Aut(G)$ such that $A^{\Inn(G)}=A$ and the hypercenter of $A/(A\cap \Inn(G))$ has finite index $k$.

\medskip\noindent
(B)\ 
If  the $A$-hypercenter of $G$ has finite index $t$, then there is a finite normal $A$-subgroup   $L$ with order  bounded by a function of\ $(t,k)$\ such that $G/L$ is $A$\nobreakdash-hypercentral.

\medskip\noindent
(H)\  If there is a finite normal $A$-subgroup $L$ with order $d$ such that $G/L$ is $A$\nobreakdash-hypercentral, then the $A$-hypercenter of $G$ has finite index  bounded by a funcion of  $(d,k)$.
\end{Teo}

Remark that this theorem generalizes Theorems 4 and 3 of \cite{DKP} where the same picture is considered, but with more restrictive conditions, that is $A$ contains $\Inn(G)$, the factor $A/\Inn(G)$ is finite and the involved series which are stabilized by $A$ are finite. Clearly, our bounding functions do not depend on the length of the considered series.  

Finally note that the hypothesis that $A$ is normalized by $\Inn(G)$ is necessary, as shown by Example in Sect. 2 below.


\section{Proof of  Theorem \ref{TeoremaCasolo}}

To prove Theorem \ref{TeoremaCasolo} we use a key lemma. Recall that we denote the hypercenter of a group $G$ by $Z_\infty(G)$.

\begin{Lemma}\label{LemmaCasolo}
Let $A\le H$ be  normal subgroups of a group $G$ with $A$  finite and $A\le Z(H)$. If $G/C_{G}(H)$ is locally nilpotent and $H/A\le Z_{\infty}(G/A)$, then $H\le Z_{\infty}(G)A$.
\end{Lemma}
\pf Arguing by induction on the order of $A$, we may assume that $A$ is minimal normal in $G$. Then $A$ is an elementary abelian $p$-group for some prime $p$. If $A\cap Z(G)\ne 1$, then $A\le Z(G)$ by minimality of $A$ and so we have $H\le Z_{\infty}(G)A$.

Suppose then $A\cap Z(G)=1$ (and so $A\cap Z_{\infty}(G)=1$) and let $N:=Z_{\infty}(G)\cap H$. Note that the hypotheses hold for the subgroups $\bar A:=AN/N$, $\bar H:=H/N$ of the group $\bar G:=G/N$. Since from $\bar H\le Z_{\infty}(\bar G)\bar A$ it follows $H\le Z_{\infty}(G)A$, we may assume $Z_{\infty}(G)\cap H=1$ .

We claim that $H=A$ (note that $H\le Z_{\infty}(G)A$ if and only if   $H=H\cap  Z_{\infty}(G)A=(H\cap  Z_{\infty}(G))A=A$). Suppose, by contradiction, $H>A$ and let $X/A\ne 1$ be either infinite cyclic or of prime order $r$ and contained in  $(H/A)\cap Z(G/A)$.  Since by hypotheses $A\le Z(H)$, then $X$ is abelian and  $X\n G$, clearly.  

 Let us show now that $X$ is a $p$-group. If, by contradiction,  $X/A$ is infinite or $r\ne p$, then $X^{p}\ne 1$ and $X^{p}\cap A=1$. Thus $X^{p}$ is $G$-isomorphic to $X^{p}A/A\le Z_{\infty}(G/A)$. Hence $X^{p}\le H\cap Z_{\infty}(G)=1$, a contradiction.  So $X/A$ has order $p$.
 
Assume, again by contradiction,  $X^{p}\ne 1$. By minimality of $A$, we have $X^{p}=A=[G,X]$ and so $[G,A]=[G,X^{p}]=[G,X]^{p}=A^{p}=1$, a contradiction.

Then $X$ is a finite  elementary abelian $p$-group. Since $[G,A]=A\le X$, the subgroup $X\rtimes (G/C_{G}(X))$ of the holomorph of $X$ is not nilpotent, and so $G/C_{G}(X)$ is not a $p$-group. Hence there are  a prime $q\ne p$ and a normal non-trivial $q$-subgroup $Q/C_{G}(X)$ of $G/C_{G}(X)$. Since $Q\not\le C_{G}(X)$,  then $[X,Q]\ne 1$. Thus  $[X,Q]=A$, as $[X,Q]\le A$ and  by minimality of $A$.

By a standard argument on coprime actions (see for example Exercise 4.1 in \cite{A}), we have 
$$X=[X,Q]\times C_{X}(Q)=A\times C_{X}(Q),$$
therefore $C_{X}(Q)\ne 1$. On the other hand,  $C_{X}(Q)$ is a normal subgroup of $G$ and so $C_{X}(Q)\le Z_{\infty}(G)\cap H=1$, a contradition which gives the claim $H=A$.
\qed

\noindent {\bf Proof of Theorem \ref{TeoremaCasolo}.}\quad
Let us apply Lemma  \ref{LemmaCasolo} with $A:=Z(L)$ and $H:=C_{G}(L)$. In fact on  one hand $H/A=H/(H\cap L)\simeq_{G} LH/L$, then $H/A\le Z_{\infty}(G/A)$. On the other hand $L\le C_{G}(H)$ and so $G/C_{G}(H)$ is hypercentral, since it is an image of $G/L$. Therefore $H\le Z_{\infty}(G)A$. Hence
$$|H/(H\cap Z_{\infty}(G))|=|A(Z_{\infty}(G)\cap H)/(Z_{\infty}(G)\cap H)|\le |A|=|Z(L)|.$$
Since $H=C_{G}(L)$, then $|G/H|\le |\Aut(L)|$. Thus
$$|G/Z_{\infty}(G)|\le |G/H|\cdot |H/(H\cap Z_{\infty}(G))|\le |\Aut(L)|\cdot |Z(L)|.$$
\vskip-7mm\qed

\noindent {\bf Proof of Corollary \ref{CorTeoremaCasolo}.}  
Note that $Z_{d+m}(G)=Z_\infty(G)$ has finite index. Thus if $d\le m$, the statement follows directly from Theorem  \ref{TeoremaCasolo}. Otherwise, $|G/Z_{2m}(G)|$ is bounded by the maximum of the $h(d,i)$ with $i=1,\ldots,d$, where $h(d,m)$ is the bounding function in Hall Theorem.
\qed

From Theorem 1 and the above quoted result from \cite{KOS} we deduce a corollary which gives a rather complete picture of finite-by-hypercentral groups.

\begin{Cor}\label{Corollario}  If $G$ is a group with  a (finite) normal series
$$G=G_0\ge F_1\ge G_1\ge\ldots\ge F_n\ge G_n=1$$\vskip-4mm 
\noindent where\\ - each factor $F_{i}/G_i$ is finite with order $t_i>1$, \\ - each factor $G_{i-1}/F_i$ is contained in the hypercenter of $G/F_i$, \\
then there is a  normal subgroup  $L$ with finite order bounded by a function of $t=t_1\cdot...\cdot t_n$ such that $G/L$ is hypercentral.

\hskip-1mm Moreover the hypercenter of $G$ has finite index bounded by a function of $t$.
\end{Cor}

\pf   Define recursively a function $f:\N\to \N$ by means of $f(1)=1$ and $f(t+1)=(t+1)g(g(f(t)))$ for each $t\in\N$, where $g(t):=t^{1+\log_2t}$.

 We show that there is $L\n G$ such that $|L|\le f(t)$ and $G/L=Z_\alpha(G/L)$ for $\alpha:=\alpha_n+\ldots+\alpha_1+m'$, where 
the  $\alpha_i$'s are ordinal numbers such that $G_{i-1}/F_i\le Z_{\alpha_i}(G/F_i)$ for each $i$ and  $m'\in\N$ may be bounded by a function of $t$ and of the $\alpha_i$'s which are finite. 
 Since $f(t)\ge t$ for each $t$, the statement is trivial if $n=1$.  

Assume then by induction on $n$ that there is a normal series $$G\ge F_{n-1}\ge G_{n-1}\ge F_n\ge G_n=1$$ such that 
$G/F_{n-1}$ is hypercentral of type $\alpha'=  \alpha_{n-1}+\ldots+\alpha_1+m''$, with $m''\in\N$ and  $|F_{n-1}/G_{n-1}|\le f(t_*)$ with $t_*=t_1\cdot...\cdot t_{n-1}$. 
Applying Theorem \ref{TeoremaCasolo} to $G/G_{n-1}$, if $Z/G_{n-1}:=Z_{(\lceil\log_2f(t_*)\rceil+\alpha')}(G/G_{n-1})$, then $|G/Z|\le g(f(t_*))$. Thus,  applying Theorem B of [KOS]  to $G/F_n$, we have that there is a normal subgroup $L$ such that $G/L$ is hypercentral with ordinal type at most $\alpha_n+\lceil \log_2f(t_*)\rceil+\alpha'+\lceil\log_2 g(f(t_*))\rceil$ and $|L/F_n|\le g(g(f(t_*)))$. We have:\ 
$|
L|\le t_1 g(g(f(t_*)))\le 
 t g(g(f(t-1)))= f(t)$, as wished. 
\QED

\noindent{\bf Remark}: In the above proof, if $\alpha$ is infinite, then clearly $G/Z_\alpha(G)$ is finite. Otherwise, if $G_{i-1}/F_i\le Z_{m_i}(G/F_i)$ for each $i$ with $m_i\in\N$, then there is a finite normal subgroup $L$ such that   $G/L=Z_{m}(G/L)$  with $m:=m_1+m_2+\ldots+m_n$, by Theorem B in \cite{FM}.
Hence, in this case,  $G/Z_{2m}(G)$ is finite.


\section{Proof of Theorem 2}

\noindent{\bf Proof of Theorem 2.} 
Let $\alpha'$ such that  $B/(A\cap \Inn(G)):=Z_{\alpha'}(A/(A\cap \Inn(G))$
 has finite index in $A/(A\cap \Inn(G))$.  Consider the subgroup $S:=G\rtimes A$ of the holomorph group of $G$.


Assume first $A\ge \Inn(G)$. Let  $G_\delta:=Z_\delta(G,A)$ for any ordinal $\delta$. {We claim:  $$(*)\hskip1cm\forall \delta\hskip0.5cm  S_\delta:=G_\delta\bar G_\delta \le Z_{\delta}(S).$$}
By induction, suppose true for $\delta$. Note that $\bar G\le A$ acts by conjugation on $G$ the same way as $G$. We have 
$[S_{\delta+1},S]=[G_{\delta+1}\bar G_{\delta+1},GA]$. On one hand, we have 
$ 
[G_{\delta+1},GA]\le[G_{\delta+1},A]\cdot[G_{\delta+1},G]^A\le G_\delta$. On the other hand,  
$ [\bar G_{\delta+1},GA]\le
[\bar G_{\delta+1},A]\cdot[ \bar G_{\delta+1},G]^{A}\le   \bar G_\delta G_\delta=S_\delta$.
It follows 
  $S_{\delta+1}\le Z_{\delta+1}(S)$ and the claim is proved since the limit ordinal step is trivial.

\medskip

To prove (B) in the case $A\ge \Inn(G)$, let $\alpha$ be such that  $Z_\alpha(G,A)$ has finite index in $G$ and  note that 
 in the normal series
$$S=GA\ge GB\ge G\bar G\ge G_\alpha\bar G_\alpha\ge 1$$
the factors $GA/ GB$ and $G\bar G/G_\alpha\bar G_\alpha$ are finite with order $k$ and $t^2$, respectively.  
Moreover, by $(*)$, factors $GB/G\bar G$ and $ G_\alpha\bar G_\alpha$ are contained in the $\alpha'$th and $\alpha$th center of $S/G\bar G$ and $S$, respectively. Thus we apply Corollary \ref{Corollario} to the group $S=GA$. Then 
the statement (for the group $G$) follows easily.   

\medskip

Concerning part (H)  in the case $A\ge \Inn(G)$, consider the normal series
$$S=GA\ge GB\ge G\bar G\ge L\bar L\ge 1.$$ Note that $GA/ GB$ and  $L\bar L$ are finite with order  $k$ and $d^2$, respectively. Moreover, if $\alpha_1$ is such that  $Z_{\alpha_1}(G/L,A)$ has finite index in $G/L$, then by $(*)$ we have that $GB/ L\bar L$ is contained in the ($\alpha_1+\alpha'$)th $A$-center of  $S/ L\bar L$. We may apply Corollary \ref{Corollario} and get the statement. 

\medskip

To deal with the more general case, let $\bar N := A\cap \Inn(G)$ such that  $Z(G)\le N\le G$. Note that  $[G,A]\le N$, as $\overline{[g,\gamma]}=  [\bar g,\gamma]\in A\cap \Inn(G)\   \forall \gamma \in A$ since $A^{\Inn(G)}=A$. Thus $A$ acts trivially on $G/N$. 
Moreover the group $\tilde A:=A/C_A(N)$ may be considered as a group of automorphisms on $N$ containing $\Inn(N)$. Thus, to prove (H), one may apply the above case to $N$ and $\tilde A:=A/C_A(N)$.



To prove (B) in the general case note that, by the above, the subgroup $Z:=Z_\infty(N,A)$ has finite index in $N$, bounded by a function of $|L\cap N|\le|L|$. Let $K/Z$ be the $A$-hypercenter of $G/Z$. Clearly,  $K\cap N=Z$. Moreover $K/Z=Z(G/Z,A)$. Consider then $C/Z:=C_{G/Z}([G,A]Z/Z)$ and note that  $C$ has finite index in $G$, since  $[G,A]\le N$.  
By applying the Three Subgroup Lemma to $A$, $C/Z$, $C/Z$, we have that $A$ acts trivially on the derived subgroup of $C/Z$. Thus $C'Z/Z\le C_{G/Z}(A)\le K/Z$. Therefore 
 $CK/K$ is abelian. We consider the series 
$$G \ge CK \ge K \ge Z \ge 1.$$ 
The index of $CK$ in $G$ is finite and bounded by a function of $d=|L|$, as $|N/Z|$ is.
Then consider the action of $A$ on the abelian group $\hat G:=CK/K$.
Since $K\cap N=Z$, we have that $|NK/K|$ is bounded by a function of $d$. Thus 
 the image of $A\cap \Inn (G)$  in $\hat A:=A/C_A(\hat G)$ is finite with order bounded by a function of $d$. By Corollary \ref{Corollario}, $Z_{\alpha'}(\hat A)$ has finite index $q$ in $\hat A$, bounded by a function of $d$ and $k$. Recall that  $\hat G$ is abelian and $[\hat G,\hat A]$ is finite, as $[G,A]$ is finite modulo $K$. Let $\hat S:=\hat G\rtimes \hat A$. Then $Z_{1+\alpha'}(\hat S/[\hat G,\hat A])$ has finite index at most $q$. By Theorem \ref{TeoremaCasolo}, the index of $Z_{1+\alpha'}(\hat S)$ in $\hat S$ is finite and bounded by a function of $d$ and $q$. Thus the $A$-hypercenter of $\hat G:=CK/K$ has finite index and bounded by a function of $d$ and $k$, as wished.\QED

\noindent{\bf Remark}: in the case  $A\ge \Inn(G)$ of the above proof, if $\alpha$,  $\alpha_1$ and $\alpha'$ are finite, we have that:\\ 
- in case (B),  the $2(\alpha+\alpha')$th $A$-center  has finite index in $G$, by the above quoted result in \cite{FM}. In particular, for $\alpha'=0$ we have Theorem 3 of \cite{DKP}.\\
- in case (H), there is a boundedly finite normal $A$-subgroup $L$ such that $G/L$ coincides with its $(\alpha_1+\alpha')$th $A$-center. This follows by applying the remarks  after Corollary \ref{Corollario}  to the group $S$. In particular, for   $\alpha'=0$ we have Theorem 2 and 4 of \cite{DKP}.

\medskip 
Let us see that the condition that $A$ is normalized by $\Inn(G)$ is necessary.

\medskip
\noindent {\bf Example} {\it There is an elementary abelian group $G$ and a bounded abelian group $A\le \Aut(G)$ such that $G/Z_{\omega}(G,A)$ is finite (of prime order), while $G/L$ is not $A$-hypercentral, for any finite  $A$-subgroup $L\le G$.}

\medskip

\pf Let $G:=Dr_{i<\omega}\la a_{i}\ra$ be an elementary abelian $p$-group, where $p$ is an odd prime and let $Z:=Dr_{0<i<\omega}\la a_{i}\ra$. For any $i>0$, consider $\g_{i}\in \Aut(G)$ centralizing $Z$, and such that $a_{0}^{\g_{i}}:=a_{0}a_{i}$. Let $\tau\in \Aut(G)$ centralizing $Z$ and such that $a_{0}^{\tau}:=a_{0}^{2}$. Let $A$ be the subgroup of $\Aut(G)$ generated by $\tau$ and all the $\g_{i}$'s. Then $Z=Z_{1}(G,A)$ has index $p$ in $G$, while if $K$ is a proper $A$-subgroup of $G$, then $a_{0}\not\in K$, as $a_{0}^{A}=G$. Clearly $\tau$ does not centralizes $a_{0}$ mod $K$. Thus $G/K$ is not $A$-hypercentral, for any proper $A$-subgroup $K$ of $G$ and in particular for any finite $A$-subgroup $L\le G$.\qed

We finish by noticing that Theorem \ref{generalizeBH} may be formulated in a different way. Recall that the factor of two consecutive terms of a series is called just factor.

\begin{Cor}\label{Corollario2}  Let $A$ be a finite-by-hypercentral group of automorphisms of a group $G$ such that $A^{\Inn(G)}=A$.

\hskip-1mm If there is an ascending normal series in $G$ with a finite number of finite factors and such that $A$ acts trivially on all other factors, then:\\ 
i)  there is 
a finite index normal $A$-subgroup $G_0$ of $G$ such that 
$A$ stabilizes an ascending $G$-series  of $G_0$;\\
ii) there is  is a  finite  normal $A$-subgroup $L$ such that 
$A$ stabilizes an ascending $G$-series of $G/L$. 
\end{Cor}


Carlo Casolo, 
Dipartimento di Matematica “U. Dini”, Universit\`a di Firenze, Viale Morgagni 67A, I-50134 Firenze, Italy. \\ email: casolo@math.unifi.it

\bigskip
Ulderico Dardano,  Dipartimento di Matematica e
Applicazioni ``R.Caccioppoli'', Universit\`a di Napoli ``Federico
II'',  Via Cintia - Monte S. Angelo, I-80126 Napoli, Italy. \\ email: dardano@unina.it

\bigskip
Silvana Rinauro,  
Dipartimento di Matematica, Informatica ed Economia, Universit\`a della
Basilicata, Via dell'Ateneo Lucano 10 - Contrada Macchia Romana,
I-85100 Potenza, Italy.\\ email: silvana.rinauro@unibas.it

\end{document}